\documentclass{amsart} 

\newtheorem{thm}{Theorem}[section]

\theoremstyle{definition}\newtheorem{question}[thm]{Question}

\theoremstyle{remark}\newtheorem{remark}[thm]{Remark}

\theoremstyle{definition}\newtheorem{example}[thm]{Example}

\def\R{{\Bbb R}}

\def\Q{{\Bbb Q}}

\def \air{{\vskip 12pt\noindent}\par}

\def \ib{\begin{enumerate}}
\def \ie{\end{enumerate}}

\def\ib {\begin{enumerate}}
\def\ie {\end{enumerate}}

\def\adots{\mathinner{\mkern2mu\raise 1pt\hbox{.}\mkern 3mu\raise
4pt\hbox{.}\mkern1mu\raise 7pt\hbox{{.}}}}

\def \Gal {{\rm Gal}}

\title{Bounding the rational sums of squares over totally real fields}
\author{Ronan Quarez}

\address{IRMAR (CNRS, URA 305), Universit\'e de Rennes 1, Campus de Beaulieu\\ 35042 Rennes Cedex, France} 
\email{e-mail : ronan.quarez@univ-rennes1.fr}
\date{\today} 


\begin{document}
\maketitle

\begin{abstract}
Let $K$ be a totally real Galois number field. C. J. Hillar proved that if $f\in\Q[x_1,\ldots,x_n]$ is a sum of $m$ squares in $K[x_1,\ldots,x_n]$, then $f$ is a sum of $N(m)$ squares in $\Q[x_1,\ldots,x_n]$, where 
$N(m)\leq 2^{[K:\Q]+1}\cdot {[K:\Q]+1 \choose 2}\cdot 4m$, the proof being constructive. \par
We show in fact that $N(m)\leq (4[K:\Q]-3)\cdot m$, the proof being constructive as well.
\end{abstract}

\section{Introduction}
In the theory of semidefinite linear programming, there is a question by Sturmfels  

\begin{question}[Sturmfels]
{If $f\in\Q[x_1,\ldots,x_n]$ is  a sum of squares in $\R[x_1,\ldots,x_n]$, then is $f$ also a sum of squares in $\Q[x_1,\ldots,x_n]$ ?}
\end{question}

Hillar (\cite{Hi}) answers the question in the case where the sum of squares has coefficients in a totally real Galois number field :

\begin{thm}[Hillar]
Let $f\in\Q[x_1,\ldots,x_n]$ be a sum of $m$ squares in $K[x_1,\ldots,x_n]$ where $K$ is a totally real Galois extension of $\Q$. Then, $f$ is a sum of  
$$2^{\displaystyle [K:\Q]+1}\cdot {[K:\Q]+1 \choose 2}\cdot 4m$$
squares in $\Q[x_1,\ldots,x_n]$.
\end{thm}

The aim of this note is to show, modifying a little bit Hillar's proof, that only $(4[K:\Q]-3)\cdot m$ squares are needed (that is Theorem \ref{betterbound}). Moreover, as in \cite{Hi}, the argument is constructive.

\section{Hillar's method}
Having in mind the Lagrange's four squares Theorem, we focus ourselves on {\it rational sum of squares} i.e. linear combination of squares with positive rationnal coefficients.\par
Let $K$ be a totally real Galois extension of $\Q$ which we write $K=\Q(\theta)$ with $\theta$ a real algebraic number, all of whose conjugates are also real. We set $r=[K:\Q]$ and $G=\Gal(K/\Q)$.\par

Let $f\in\Q[x_1,\ldots,x_n]$ be a sum of $m$ squares in $K[x_1,\ldots,x_n]$, namely $f=\sum_{k=1}^mf_k^2$, with $f_k\in K[x_1,\ldots,x_n]$.
Summing over all actions of $G$ (i.e. ``averaging"), we get 

\begin{equation}\label{averaging}
f=\frac{1}{\mid G\mid}\sum_{k=1}^m\sum_{\sigma\in G}(\sigma f_k)^2
\end{equation}

Next, we write each $f_k$ in the form 
$$f_k=\sum_{i=0}^{r-1}q_i\theta^i$$
where $q_i\in\Q[x_1,\ldots,x_n]$. Then, 
\begin{equation}\label{decomp}
\sum_{\sigma\in G}(\sigma f_k)^2=\sum_{j=1}^r\left(\sum_{i=0}^{r-1}q_i(\sigma_j\theta)^i\right)^2
\end{equation}

We may write this sum of squares as the following product of matrices 

$$\left(\begin{array}{c}
q_0\\
\vdots\\
q_{r-1}\\
\end{array}\right)^T\left(\begin{array}{cccc}
1&\sigma_1\theta&\ldots&(\sigma_1\theta)^{r-1}\\
\vdots&\vdots&\ddots&\vdots\\
1&\sigma_r\theta&\ldots&(\sigma_r\theta)^{r-1}\\
\end{array}\right)^T
\left(\begin{array}{cccc}
1&\sigma_1\theta&\ldots&(\sigma_1\theta)^{r-1}\\
\vdots&\vdots&\ddots&\vdots\\
1&\sigma_r\theta&\ldots&(\sigma_r\theta)^{r-1}\\
\end{array}\right)\left(\begin{array}{c}
q_0\\
\vdots\\
q_{r-1}\\
\end{array}\right)$$

We obtain what is called a Gram matrix (cf \cite{CLR}) associated to the sum of squares in (\ref{decomp}). Let

$$G=\left(\begin{array}{cccc}
1&\sigma_1\theta&\ldots&(\sigma_1\theta)^{r-1}\\
\vdots&\vdots&\ddots&\vdots\\
1&\sigma_r\theta&\ldots&(\sigma_r\theta)^{r-1}\\
\end{array}\right)^T
\left(\begin{array}{cccc}
1&\sigma_1\theta&\ldots&(\sigma_1\theta)^{r-1}\\
\vdots&\vdots&\ddots&\vdots\\
1&\sigma_r\theta&\ldots&(\sigma_r\theta)^{r-1}\\
\end{array}\right)$$
Note that the entries of $G$ are in $\Q$ since they are invariant under the $\sigma_j$'s. \air

Now, we come to the slight modification of the proof of Hillar that will improve the bound.

%%%%%%%%%%
\section{LU-decomposition of the Gram matrix}
If $u(x)$ denotes the minimal polynomial of the Galois extension $K$ over $\Q$, then the $(i,j)$-th entry of the matrix $G$ is the $i+j-2$-th Newton sum of $(\sigma_1,\ldots,\sigma_r)$ the roots of $u(x)$. It is well known that the rank of $G$ is equal to $r$ and its signature (the difference between the positive eigenvalues and the negative ones) is equal to the number of real roots of $u(x)$ (see for instance \cite[Theorem 4.57]{BPR}).
In our case, we readily deduce that $G$ is a positive definite matrix since $K$ is totally real. Thus, all its principal minors are different from zero (they are strictly positive !) and $G$ admits a LU-decomposition which we may put in a symmetric form 

$$G=U^TDU$$
where $D$ is diagonal and $U$ is upper triangular with diagonal identity, and $U,D$ have rational entries. \par
We may view this decomposition as a matricial realization of the Gauss algorithm which reduce the quadratic form given by $G$.
\air
Now, if we denote by $f_1,\ldots,f_r$ the polynomials in $\Q[x_1,\ldots,x_r]$ appearing as the rows of the matrix  $U\times\left(\begin{array}{c}
q_0\\
\vdots\\
q_{r-1}\\
\end{array}\right)$ and by $d_1,\ldots,d_r$ the rational entries onto the diagonal of $D$, then we get from (\ref{decomp}) the identity :

\begin{equation}\label{afterLU}
\frac{1}{\mid G\mid}\sum_{\sigma\in G}(\sigma f_k)^2=\frac{d_1}{\mid G\mid}g_1^2+\ldots+\frac{d_r}{\mid G\mid}g_r^2
\end{equation} 

This construction leads to 

\begin{thm}\label{betterbound}
Let $f\in\Q[x_1,\ldots,x_n]$ be a sum of $m$ squares in $K[x_1,\ldots,x_n]$, where $K$ is a totally real Galois extension of $\Q$. Then, $f$ is a sum of  
$(4[K:\Q]-3)\cdot m$
squares in $\Q[x_1,\ldots,x_n]$.
\end{thm}
 
\begin{proof} 
By (\ref{averaging}) and (\ref{afterLU}), it suffices to apply Lagrange's four squares Theorem to get that $f$ is a sum of $4[K:\Q]\cdot m$ squares in $\Q[x_1,\ldots,x_n]$.\par
But let us note that the first diagonal entry of $D$ is always $d_1=r=[K:\Q]$. Then, by the averaging process the first coefficient appearing in the rational sum of squares in (\ref{afterLU}) is $\frac{d_1}{\mid G\mid}=1$ : already a square in $\Q$ ! Whereas the others coefficients $\frac{d_i}{\mid G\mid}$ in the rational sum of squares could be any positive rational which we rewrite as a sum of $4$ squares. This concludes the proof.
\end{proof}

\remark{Beware that if we perform the Cholesky algorithm to the matrix $G$ instead of the LU-decomposition, it yields a factorisation $G=U^TU$ where $U$ is lower triangular but with entries in $\Q[\sqrt{d_1},\ldots,\sqrt{d_r}]$ for some integers $d_1,\ldots,d_r$. Then, an averaging argument would produce identities over $\Q$ but will raise the number of squares by an unexpected multiplicative factor $2^{[K:\Q]}$.} \endremark

Let us consider as an example, the simple case of quadratic extensions :
\example{
Let $K=\Q(\sqrt{d})$ where $d\in\Q$ is not a square. The extension $K$ is always Galois, and it is totally real if $d\geq0$.\par
Let $f\in\Q(x_1,\ldots,x_n)$ be such that $f=\sum_{k=1}^m(a_k+b_k\sqrt{d})^2$ with $a_k,b_k\in\Q(x_1,\ldots,x_n)$. Since $f$ has rational coefficients, by averaging we get $$f=\frac{1}{2}\sum_{k=1}^m(a_k+b_k\sqrt{d})^2+(a_k-b_k\sqrt{d})^2=\sum_{k=1}^m(a_k^2+d b_k^2)$$
It remains to write $d$ as a sum of $l\leq 4$ squares of rationals, and we get that $f$ is a sum of at most $(1+l)\cdot m$ squares in $\Q(x_1,\ldots,x_n)$.
}\endexample

As another illustration, we apply our method to \cite[Example 1.7]{Hi}  :
\example{Consider the polynomial
$$f=3-12y-6x^3+18y^2+3x^6+12x^3y-6xy^3+6x^2y^4$$ which is the following sum of squares
$$f=(x^3+\alpha^2y+\beta xy^2-1)^2+(x^3+\beta^2y+\gamma xy^2-1)^2+(x^3+\gamma^2y+\alpha xy^2-1)^2$$
in $\Q(\alpha)[x,y]$ where $\alpha,\beta,\gamma$ are the real roots of the polynomial $u(x)=x^3-3x+1$.\par
We do not need to average and directly compute the matrix $G$ and its symmetric LU-decompositon 
$$\left(\begin{array}{crr}
3&0&6\\
0&6&-3\\
6&-3&18
\end{array}\right)=\left(\begin{array}{crr}
1&0&0\\
0&1&0\\
2&-\frac{1}{2}&1
\end{array}\right)\left(\begin{array}{crr}
3&0&0\\
0&6&0\\
0&0&\frac{9}{2}
\end{array}\right)\left(\begin{array}{crr}
1&0&2\\
0&1&-\frac{1}{2}\\
0&0&1
\end{array}\right)$$
Because of the relations $\beta=2-\alpha-\alpha^2$ and $\gamma=\alpha^2-2$, the vector of polynomials $q=(q_0,q_1,q_2)^T$ is $q=(x^3+2xy^2-1,-xy^2,y-xy^2)^T$ and hence

$$f=3\left((x^3+2xy^2-1)+2\frac{}{}(y-xy^2)\right)^2+6\left(-xy^2-\frac{1}{2}(y-xy^2)\right)^2+\frac{9}{2}\left(\frac{}{}y-xy^2\right)^2$$
 
a rationnal sum of $3$ squares, to compare with the rational sum of $6$ squares obtained in  \cite{Hi}. 

}
\endexample

%%%%%%%%%%%%%%%%%%%%%%%%%%%%%%
%%%%%%%%%%%%%%%%%%%%%%%%%%%%%%
%%%%%%%%%%%%%%%%%%%%%%%%%%%%%%

\end{document}